\documentclass[reqno]{amsart}
\usepackage{amsmath}
\usepackage{xspace}
\usepackage{alltt}

\theoremstyle{definition}

\theoremstyle{plain}
\newtheorem{thm}{Theorem}
\newtheorem{crl}{Corollary}
\newtheorem{lem}{Lemma}

\numberwithin{equation}{section}
\allowdisplaybreaks[1]

\newcommand{\qBin}[3]{\genfrac{[}{]}{0pt}{0}{#1}{#2}_{#3}}

\newcommand{\qhg}{$q$\nobreakdash-hyper\-geometric\xspace}

\newlength{\axellength}

\def\eqn#1{\eqref{eq:#1}}

\begin{document}


\makeatletter
\def\print@eqnum{\tagform@{\normalsize\theequation}}
\def\endmathdisplay@a{%
  \if@eqnsw \gdef\df@tag{\tagform@{\normalsize\theequation}}\fi
  \if@fleqn \@xp\endmathdisplay@fleqn
  \else \ifx\df@tag\@empty \else \veqno \alt@tag \df@tag \fi
    \ifx\df@label\@empty \else \@xp\ltx@label\@xp{\df@label}\fi
  \fi
  \ifnum\dspbrk@lvl>\m@ne
    \postdisplaypenalty -\@getpen\dspbrk@lvl
    \global\dspbrk@lvl\m@ne
  \fi
}
\makeatother

\title[\tiny{New polynomial analogues of Jacobi's triple product and Lebesgue's identities}]
      {New polynomial analogues of Jacobi's \\
       triple product and Lebesgue's identities}
       
\author[K.~Alladi and A.~Berkovich]{Krishnaswami Alladi and Alexander Berkovich}
\address{Department of Mathematics, The University of Florida,
         Gainesville, FL~32611, USA}
\email{alladi@math.ufl.edu}
\email{alexb@math.ufl.edu}
\thanks{Research supported in part by National Science Foundation Grant DMS--0088975}

\begin{abstract}

In a recent paper by the authors, a bounded version of G\"ollnitz's (big) partition theorem 
was established. Here we show among other things how this theorem leads to nontrivial new 
polynomial analogues of certain fundamental identities of Jacobi and Lebesgue. We also derive 
a two parameter extension of Jacobi's famous triple product identity.

\end{abstract}

\maketitle

\section{Introduction}

\medskip

Jacobi's celebrated triple product identity 
\begin{equation} 
\sum_{j=-\infty}^\infty z^jq^{\frac{j(j-1)}{2}}=
\prod_{i=1}^\infty(1+zq^{i-1})(1+z^{-1}q^i)(1-q^i)
\label{eq:1.1}
\end{equation}
is one of the most fundamental results in the theory of theta functions. Employing standard 
\qhg notations
\begin{equation} 
(a_1,a_2,\ldots,a_r;q)_N=(a_1,a_2,\ldots,a_r)_N=(a_1;q)_N (a_2;q)_N (a_r;q)_N,
\label{eq:1.2}
\end{equation}
where
\begin{equation} 
(a;q)_N=(a)_N=
\begin{cases} \prod_{j=0}^{N-1}(1-aq^j), &\mbox{ if }N>0, \\ 1, &\mbox{ if }N=0, 
\end{cases}  
\label{eq:1.3}
\end{equation}
and
\begin{equation} 
(a)_\infty=\lim_{N\rightarrow\infty}(a)_N, \mbox{ for } |q|<1,
\label{eq:1.4}
\end{equation}
we can rewrite \eqn{1.1} in the equivalent form as
\begin{equation} 
\sum_{j\ge 0}q^{T_j}\frac{z^{-j}+z^{1+j}}{1+z}=(-qz,-qz^{-1},q;q)_\infty,
\label{eq:1.5}
\end{equation}
with
\begin{equation} 
T_j=\frac{j(j+1)}{2}.
\label{eq:1.6}
\end{equation}

Performing changes $q\rightarrow q^2$ and $z\rightarrow zq$ in \eqn{1.1}, it is easy to obtain
\begin{equation} 
\sum_{j=-\infty}^\infty z^jq^{j^2}=(q^2,-qz,-qz^{-1};q^2)_\infty.
\label{eq:1.7}
\end{equation}
A polynomial generalization of \eqn{1.7} due to MacMahon in \cite{M} is
\begin{equation} 
\sum_{j=-\infty}^\infty z^jq^{j^2}\qBin{2L}{L+j}{q^2}=(-qz,qz^{-1};q^2)_L,
\label{eq:1.8}
\end{equation}
where the $q$-binomial coefficients are defined as
\begin{equation} 
\qBin{n+m}{n}{q}= 
\begin{cases} \frac{(q)_{n+m}}{(q)_n(q)_m}, &\mbox{ if }n,m\ge 0, \\ 0, &\mbox{ otherwise}. 
\end{cases}  
\label{eq:1.9}
\end{equation}
Thanks to
\begin{equation} 
\lim_{L\rightarrow\infty}\qBin{2L}{L+j}{q^2}= 
\frac{1}{(q^2;q^2)_\infty},
\label{eq:1.10}
\end{equation}
\eqn{1.8} reduces to \eqn{1.7} as $L$ tends to infinity.

Andrews \cite{A1} pointed out that the proof of \eqn{1.8} requires only Euler's identity
\begin{equation} 
\sum_{j=0}^L q^{T_j}z^j\qBin{L}{j}{q}=(-qz;q)_L.
\label{eq:1.11}
\end{equation}
In this sense, \eqn{1.8} is a very simple identity. In contrast, in \cite{AB} we stated a 
rather nontrivial polynomial analogue of \eqn{1.5}, namely
\begin{equation} 
\sum_{n=0}^L q^{T_n}\frac{z^{-n}+z^{1+n}}{1+z}=
\sum_{i,j,k\ge 0} q^{T_i+T_j+T_k}z^{i-j}(-1)^k
\qBin{L-i}{j}{q} \qBin{L-j}{k}{q} \qBin{L-k}{i}{q}.  
\label{eq:1.12}
\end{equation}

In \cite{AB} we indicated that \eqn{1.12} can be deduced from a special case of the bounded 
version of G\"ollnitz's (big) partition theorem. However, no details were given there. 
Berkovich and Riese \cite{BR} used the computer algebra package \texttt{qMultiSum} \cite{R} 
to prove \eqn{1.12} by showing that both sides there satisfy identical recurrences of order 4,
namely
$$
zS_{L+4}(z)+(-z-q^{L+4}+zq^{L+4}-z^2q^{L+4})S_{L+3}(z)+(1-z+z^2)q^{L+4}(1-q^{L+3})
$$
\begin{equation} 
S_{L+2}(z)+q^{2L+7}(1-z+z^2+zq^{L+2})S_{L+1}(z)-zq^{3L+9}S_L(z)=0,
\label{eq:1.13}
\end{equation}
and by verifying initial conditions for $L\in \{0,1,2,3\}$.

One of our primary goals here is to understand precisely the \qhg status 
of \eqn{1.12}. In fact, in the next section we will show that \eqn{1.12} is a consequence 
of the somewhat mysterious Sears--Carlitz transformation of a terminating well-poised 
$_3\phi_2$-series. In \S 3 we will briefly review the refinement of the G\"ollnitz's theorem 
due to Alladi, Andrews and Gordon \cite{AAG} and its bounded version recently discovered 
by Alladi and Berkovich \cite{AB}. In \S 4 we will reinterpret the generalized G\"ollnitz 
theorem in \cite{AAG} and \cite{AB} as a weighted partition identity. 
As a fallout of this reinterpretation we will immediately obtain \eqn{1.12} in \S 4, and then 
derive in \S 5 a two-parameter extension of \eqn{1.1} in a form of a double series. 

One of the fundamental results in the theory of partitions and $q$-series is Lebesgue's 
identity \cite{A1}, which we state in the following form:
\begin{equation}
   \sum_{i\ge 0}q^{T_i}\frac{(z^2 q)_i}{(q)_i}=(z^2q^2;q^2)_\infty(-q)_\infty.
\label{eq:4.20}
\end{equation}
Section 6 is devoted to the proof of the following new polynomial analogue of \eqn{4.20}:
$$
\sum_{i,j\ge 0} q^{T_i+T_j}(-z^2)^j \qBin{L-j}{i}{q} \qBin{i}{j}{q}=
$$
\begin{equation} 
\sum_{i,j,k\ge 0} (-1)^jz^{i+j}q^{T_i+T_j+T_k}
\qBin{L-i}{j}{q} \qBin{L-j}{k}{q} \qBin{L-k}{i}{q}.  
\label{eq:1.14}
\end{equation}

In \S 7 we briefly describe prospects for future work. Finally, some technical details 
pertaining to the determination of the weights in Theorem~4 can be found in the Appendix.

We conclude this section by recalling some standard \qhg definitions along with 
selected summation and transformation formulas \cite{GR}.

The generalized basic hypergeometric function $_{r+1}\phi_r$ is defined as
\begin{equation} 
{}_{r+1}\phi_r \left(\begin{array}{l} a_1,a_2,\ldots,a_{r+1};q,z \\
b_1,b_2,\ldots,b_r \end{array}\right)=
\sum_{j=0}^\infty\frac{(a_1,a_2,\ldots,a_{r+1})_j}{(b_1,b_2,\ldots,b_r,q)_j}z^j.
\label{eq:1.15}
\end{equation}
The $q$-binomial theorem (Cauchy's identity) states that
\begin{equation} 
{}_1\phi_0 \left(\begin{array}{c} a;q,z \\ - \end{array}\right)=
\sum_{n\ge 0}\frac{(a)_n}{(q)_n}z^n=\frac{(az)_\infty}{(z)_\infty}.
\label{eq:1.16}
\end{equation}
The identities  
\begin{equation} 
{}_2\phi_1 \left(\begin{array}{l} a,q^{-n};q,q \\ c \end{array}\right)=
a^n\frac{(\frac{c}{a})}{(c)_n},
\label{eq:1.17}
\end{equation}
and
\begin{equation} 
{}_2\phi_1 \left(\begin{array}{l} a,b;q,-\frac{q}{b} \\ \frac{aq}{b} \end{array}\right)=
(-q)_\infty 
\frac{(aq,\frac{aq^2}{b^2};q^2)_\infty}{(-\frac{q}{b},\frac{aq}{b};q)_\infty}
\label{eq:1.18}
\end{equation}
are known as $q$-Chu--Vandermonde sum and $q$-Kummer (Bailey--Daum) sum, respectively. 
Heine's third transformation for $_2\phi_1$ can be written as
\begin{equation} 
{}_2\phi_1 \left(\begin{array}{l} a,b;q,z \\ c \end{array}\right)=
\frac{(\frac{abz}{c})_\infty}{(z)_\infty}
{}_2\phi_1 \left(\begin{array}{l} \frac{c}{a},\frac{c}{b};q,\frac{abz}{c} \\ 
c \end{array}\right).
\label{eq:1.19}
\end{equation}
Finally, the Sears--Carlitz transformation of a terminating $_3\phi_2$ series is
\begin{equation} 
{}_3\phi_2 
\left(\begin{array}{l} a,b,c;q,\frac{aqz}{bc} \\ \frac{aq}{b},\frac{aq}{c} \end{array}\right)=
\frac{(az)_\infty}{(z)_\infty}
{}_5\phi_4 
\left(\begin{array}{l}\sqrt{a},-\sqrt{a},\sqrt{aq},-\sqrt{aq},\frac{aq}{bc};q,q \\ 
\frac{aq}{b},\frac{aq}{c},az,\frac{q}{z} \end{array}\right),
\label{eq:1.20}
\end{equation}
provided that $a=q^{-n}$ and $n$ is a nonnegative integer.

\medskip

\section{$q$-hypergeometric approach to \eqn{1.12}}

\medskip

We begin by rewriting the right hand side (RHS) of \eqn{1.12} as
$$
\sum_{i,j,k\ge 0}q^{T_i+T_j+T_k}z^{i-j}(-1)^k\qBin{L-i}{j}{q}\qBin{L-j}{k}{q}\qBin{L-k}{i}{q}=
$$
\begin{equation} 
\sum_{\substack{i,j\ge 0 \\ i+j\le L}}q^{T_i+T_j}z^{i-j}\frac{(q)_L}{(q)_i(q)_j(q)_{L-i-j}}
{}_2\phi_1
\left(\begin{array}{l}q^{i-L},q^{j-L};q,q^{1+L-i-j} \\ q^{-L} \end{array}\right).
\label{eq:2.1}
\end{equation}
Remarkably, $_2\phi_1$ with the choices of parameters given above can be evaluated to be 
a product. To see this we first make use of \eqn{1.19} to transform $_2\phi_1$ as
\begin{equation} 
{}_2\phi_1
\left(\begin{array}{l}q^{i-L},q^{j-L};q,q^{1+L-i-j} \\ q^{-L} \end{array}\right)=
(q)_{L-i-j}{}_2\phi_1
\left(\begin{array}{l}q^{-i},q^{-j};q,q \\ q^{-L} \end{array}\right),
\label{eq:2.2}
\end{equation}
and then employ \eqn{1.17} to obtain
\begin{equation} 
{}_2\phi_1
\left(\begin{array}{l}q^{i-L},q^{j-L};q,q^{1+L-i-j} \\ q^{-L} \end{array}\right)=
(q)_{L-i-j}\frac{(q^{i-L})_j}{(q^{-L})_j},
\label{eq:2.3}
\end{equation}
where it is understood that $0\le i\le L$, $0\le j\le L$, $i+j\le L$ in \eqn{2.2} and \eqn{2.3}.
So, the RHS of \eqn{1.12} becomes
\begin{equation} 
\sum_{\substack{i,j\ge 0 \\ i+j\le L}}q^{T_i+T_j}z^{i-j}
\frac{(q)_L(q^{i-L})_j}{(q)_i(q)_j(q^{-L})_j}.
\label{eq:2.4}
\end{equation}
Now, since
\begin{equation} 
T_i+T_{j}-ij=T_{i-j}+j,
\label{eq:2.5}
\end{equation}
we can perform a change of the summation variables in \eqn{2.4} to put it in the form
\begin{equation} 
\mbox{RHS }\eqn{1.12}=\sum_{\substack{-L\le i\le L \\ 0\le j\le L}}q^{T_i+j}z^i
\frac{(q)_L(q^{i+j-L})_j}{(q)_{i+j}(q)_j(q^{-L})_j}.
\label{eq:2.6}
\end{equation}

Next, we split the sum in \eqn{2.6} as
$$
\sum_{\substack{-L\le i\le L \\ 0\le j\le L}}=\sum_{\substack{0\le i\le L \\ 0\le j\le L}}+
\sum_{\substack{-L\le i\le 0 \\ 0\le j\le L}}-\sum_{\substack{i=0 \\ 0\le j\le L}}.
$$
Performing the change $i\rightarrow -i$ in the second sum we get
$$
\mbox{RHS }\eqn{1.12}=\sum_{\substack{0\le i\le L \\ 0\le j\le L}}q^{T_i+j}z^i
\frac{(q)_L(q^{i+j-L})_j}{(q)_{i+j}(q)_j(q^{-L})_j}+
$$
\begin{equation} 
\sum_{\substack{0\le i\le L \\ 0\le j\le L}}q^{T_i+j}z^{-i}
\frac{(q)_L(q^{i+j-L})_j}{(q)_{i+j}(q)_j(q^{-L})_j}-\sum_{j=0}^L q^j 
\frac{(q)_L(q^{j-L})_j}{(q)_j^2(q^{-L})_j}.
\label{eq:2.7}
\end{equation}

We deal with the last sum in \eqn{2.7} first. To this end we use
\begin{equation} 
(q^{j-L})_j=\frac{(q^{-L})_{2j}}{(q^{-L})_j}=
\frac{(q^{-\frac{L}{2}},-q^{-\frac{L}{2}},q^{\frac{1-L}{2}},-q^{\frac{1-L}{2}})_j}{(q^{-L})_j},
\label{eq:2.8}
\end{equation}
with $j\le L$ to derive
$$
\sum_{j=0}^L q^j \frac{(q)_L(q^{j-L})_j}{(q)_j^2(q^{-L})_j}=
(q)_L\lim_{\substack{z\rightarrow 1 \\ b\rightarrow q}}{}_4\phi_3
\left(\begin{array}{l} 
q^{-\frac{L}{2}},-q^{-\frac{L}{2}},q^{\frac{1-L}{2}},-q^{\frac{1-L}{2}};q,q
\\ q^{-L}\frac{q}{b},q^{-L}z,\frac{q}{z} \end{array}\right)=
$$
\begin{equation} 
(q)_L \lim_{\substack{c\rightarrow\infty \\ z\rightarrow 1 \\ b\rightarrow q}}\mbox{}_5\phi_4
\left(\begin{array}{l} 
q^{-\frac{L}{2}},-q^{-\frac{L}{2}},q^{\frac{1-L}{2}},
-q^{\frac{1-L}{2}},\frac{q^{1-L}}{bc};q,q \\
q^{-L}\frac{q}{b},q^{-L}\frac{q}{c},q^{-L}z,\frac{q}{z} \end{array}\right).
\label{eq:2.9}
\end{equation}
Fortunately, $_5\phi_4$ in the above equation can be transformed with the aid of 
Sears--Carlitz formula \eqn{1.20} as 
\begin{equation} 
\lim_{\substack{c\rightarrow\infty \\ z\rightarrow 1 \\ b\rightarrow q}}\mbox{}_5\phi_4=
\lim_{\substack{c\rightarrow\infty \\ z\rightarrow 1 \\ b\rightarrow q}}
\frac{(z)_\infty}{(zq^{-L})_\infty}\mbox{}_3\phi_2
\left(\begin{array}{l} 
q^{-L},b,c;q,\frac{q^{1-L}z}{bc} \\ q^{-L}\frac{q}{b},q^{-L}\frac{q}{c} \end{array}\right)=
\frac{1}{(q^{-L})_L}\sum_{n=0}^L (-1)^n q^{T_{n-1}-Ln}.
\label{eq:2.10}
\end{equation}
We remark that the limiting procedure in \eqn{2.10} is quite delicate. 
In particular, observe that
\begin{equation} 
\lim_{b\rightarrow q}\frac{(q^{-L})_n}{(\frac{q}{b}q^{-L})_n}=
\begin{cases} 1, &\mbox{ if } 0\le n\le L, \\ 0, &\mbox{ if } n>L. 
\end{cases}  
\label{eq:2.11}
\end{equation}
Now equations \eqn{2.9}, \eqn{2.10} along with
\begin{equation} 
(q^{-L})_L=(q)_L(-1)^Lq^{-T_L} 
\label{eq:2.12}
\end{equation}
imply that
$$
\sum_{j=0}^L q^j \frac{(q)_L(q^{j-L})_j}{(q)_j^2(q^{-L})_j}=
\sum_{n=0}^L (-1)^{L-n} q^{T_L+T_{n-1}-Ln}=
$$
\begin{equation} 
\sum_{n=0}^L (-1)^{L-n} q^{T_{L-n}}=\sum_{n=0}^L (-1)^n q^{T_n}.
\label{eq:2.13}
\end{equation}
That concludes our treatment of the last sum in \eqn{2.7}.

We now turn our attention to the first sum in \eqn{2.7}. We start by rewriting it in the form
\begin{equation} 
(q)_L \sum_{i=0}^L \frac{z^i q^{T_i}}{(q)_i}
\lim_{\substack{z\rightarrow q^{-i} \\ b\rightarrow q}}\mbox{}_4\phi_3
\left(\begin{array}{l} 
q^{\frac{i-L}{2}},-q^{\frac{i-L}{2}},q^{\frac{1+i-L}{2}},-q^{\frac{1+i-L}{2}};q,q \\
q^{i-L}\frac{q}{b},q^{i-L}z,\frac{q}{z} \end{array}\right).
\label{eq:2.14}
\end{equation}
As before, we will treat $_4\phi_3$ in \eqn{2.14} as a limiting case of $_5\phi_4$
\begin{equation} 
\lim_{\substack{z\rightarrow q^{-i} \\ b\rightarrow q}}\mbox{}_4\phi_3=
\lim_{\substack{c\rightarrow\infty \\ z\rightarrow q^{-i} \\ b\rightarrow q}}\mbox{}_5\phi_4
\left(\begin{array}{l} 
q^{\frac{i-L}{2}},-q^{\frac{i-L}{2}},q^{\frac{1+i-L}{2}},
-q^{\frac{1+i-L}{2}},\frac{q^{1+i-L}}{bc};q,q \\
q^{i-L}\frac{q}{b},q^{i-L}\frac{q}{c},q^{i-L}z,\frac{q}{z} \end{array}\right),
\label{eq:2.15}
\end{equation}
to which the Sears--Carlitz formula \eqn{1.20} can be applied. This way after simplification 
we have for the first sum in \eqn{2.7}
$$
\sum_{i\ge 0} z^i q^{T_L}(-1)^{L-i} 
\sum_{n=0}^{L-i}(-1)^n q^{T_{n-1}-Ln}=
\sum_{n=0}^L (-1)^{L-n} q^{T_{L-n}}
\sum_{i=0}^{L-n}(-z)^i=
$$
\begin{equation} 
\sum_{n=0}^L (-1)^n q^{T_n}\frac{1+(-1)^nz^{n+1}}{1+z}.
\label{eq:2.16}
\end{equation}
Thus, we have established that
\begin{equation} 
\sum_{L\ge i,j\ge 0} z^i q^{T_i+j}\frac{(q)_L(q^{i+j-L})_j}{(q)_{i+j}(q)_j(q^{-L})_j}= 
\sum_{n=0}^L (-1)^n q^{T_n}\frac{1+(-1)^nz^{n+1}}{1+z}.
\label{eq:2.17}
\end{equation}
We now note that the penultimate sum in \eqn{2.7} can be easily obtained from the 
first one there by $z\rightarrow\frac{1}{z}$. Hence,
\begin{equation} 
\sum_{0\le i,j\le L} z^{-i} q^{T_i+j}\frac{(q)_L(q^{i+j-L})_j}{(q)_{i+j}(q)_j(q^{-L})_j}= 
\sum_{n=0}^L (-1)^n q^{T_n}\frac{z+(-1)^n z^{-n}}{1+z}.
\label{eq:2.18}
\end{equation}
Finally, combining \eqn{2.7}, \eqn{2.13}, \eqn{2.17} and \eqn{2.18} we see that
$$
\mbox{RHS }\eqn{1.12}=\sum_{n=0}^L (-1)^n q^{T_n}
\left\{\frac{1+(-1)^n z^{n+1}}{1+z}+\frac{z+(-1)^n z^{-n}}{1+z}-\frac{1+z}{1+z}\right\}=
$$
\begin{equation} 
\sum_{n=0}^L q^{T_n}\frac{z^{-n}+z^{n+1}}{1+z}=\mbox{ LHS }\eqn{1.12},
\label{eq:2.19}
\end{equation}
as desired.

It is to be noted that when $L\rightarrow\infty$ the right hand side of \eqn{1.12} becomes
\begin{equation} 
\sum_{i,j,k\ge 0} \frac{q^{T_i+T_j+T_k}}{(q)_i(q)_j(q)_k}z^{i-j}(-1)^k=
(-qz,-qz^{-1},q;q)_\infty,
\label{eq:2.20}
\end{equation}
thanks to a limiting case of \eqn{1.11}:
\begin{equation} 
\sum_{i\ge 0} \frac{q^{T_i}}{(q)_i}z^i=(-qz)_\infty.
\label{eq:2.21}
\end{equation}

\medskip

\section{G\"ollnitz's Partition Theorem and its refinements}

\medskip

In 1967, G\"ollnitz~\cite{G} proved the following deep result:

\begin{thm} (G\"ollnitz)\\
Let $A(N)$ denote the number of partitions of $N$ in the form $N=n_1+n_2+n_3+\dots$,
such that no part is $=$ $1$ or $3$, and $n_i-n_{i+1}\ge 6$ with strict inequality if
$n_i\equiv 0,1,3 \pmod 6$.\\
Let $B(N)$ denote the number of partitions of $N$ into distinct parts $\equiv 2,4,5 \pmod 6$.
Then
$$
   A(N)=B(N).
$$
\label{thm:1}
\end{thm}

Alladi, Andrews, and Gordon \cite{AAG} reformulated and refined the above theorem using
the language of colored integers.
To state their results, we will need a few definitions.

Let $P(N,i,j,k)$ denote the number of partitions of $N$ into parts occurring in three
(primary) colors ordered as
\begin{equation} 
   \mathbf{A}<\mathbf{B}<\mathbf{C},
\label{eq:3.1}
\end{equation}
such that parts in the same color are distinct and the number of parts in colors
$\mathbf{A},\mathbf{B},\mathbf{C}$ is equal to $i,j,k$, respectively.

Next, consider partitions $\pi$, such that parts $=$ $1$ may occur in three primary
colors, namely \eqn{3.1}, but parts $\ge 2$ may occur in six colors: the three primary colors, 
and the three secondary colors $\mathbf{AB},\mathbf{AC},\mathbf{BC}$ ordered as
\begin{equation} 
	\mathbf{AB}<\mathbf{AC}<\mathbf{A}<\mathbf{BC}<\mathbf{B}<\mathbf{C}.
\label{eq:3.2}
\end{equation}
We also require that the gap between adjacent parts of $\pi$ be $\ge 1$, where the gap
may equal $1$ only if both parts are either of the same primary color or the larger part 
is in a color of higher order according to \eqn{3.2}.
These partitions $\pi$ were termed Type--$1$ partitions in \cite{AAG}.
We can now state
\begin{thm} (Alladi, Andrews, Gordon)\\
Let $G(N,a,b,c,ab,ac,bc)$ denote the number of Type--$1$ partitions
of $N$ with $a$ parts in color $\mathbf{A},\ldots,bc$ parts in color
$\mathbf{BC}$. (Note that $bc$ is a parameter and is not equal to $b$ times $c$ with similar 
interpretations for $ab$ and $ac$).
Then
\begin{equation}
\sum_{\substack{i,j,k \\ \text{constraints}}} G(N,a,b,c,ab,ac,bc)=P(N,i,j,k),
\label{eq:3.3}
\end{equation}
where the sum on the left is over the variables $a,b,c,ab,ac,bc$ subject to the
``$i,j,k$-constraints", which here and everywhere are
\begin{equation}
   \begin{cases} i=a+ab+ac,\\ j=b+ab+bc,\\ k=c+ac+bc. \end{cases}
\label{eq:3.4}
\end{equation}
\label{thm:2}
\end{thm}
To see the relation between Theorem~\ref{thm:1} and Theorem~\ref{thm:2},
we denote the integer $n$ in color $\mathbf{A}$
as $\mathbf{A}_n,\ldots$ and the integer $n$ in color $\mathbf{BC}$ as $\mathbf{BC}_n$.
Next, we replace the colored integers by the regular integers as follows:
\begin{equation}
   \left\{ \hspace*{-\arraycolsep} \begin{array}{rl} 
   \mathbf{A}_n \rightarrow 6n-4, & n\ge 1, \\
   \mathbf{B}_n \rightarrow 6n-2, & n\ge 1, \\
   \mathbf{C}_n \rightarrow 6n-1, & n\ge 1, \\[1ex]
   \mathbf{AB}_n \rightarrow 6n-6, & n>1,\\
   \mathbf{AC}_n \rightarrow 6n-5, & n>1,\\
   \mathbf{BC}_n \rightarrow 6n-3, & n>1.\end{array}\right.
\label{eq:3.5}
\end{equation}
Observe that under replacement \eqn{3.5} the ordering
$$
\mathbf{AB}_n<\mathbf{AC}_n<\mathbf{A}_n<\mathbf{BC}_n<\mathbf{B}_n<\mathbf{C}_n,\quad n>1
$$
becomes the conventional ordering 
$$
6n-6<6n-5<6n-4<6n-3<6n-2<6n-1, \quad n>1.
$$
Since part $1$ can occur only as $\mathbf{A}_1,\mathbf{B}_1,\mathbf{C}_1$, we conclude
that no conventional (uncolored) part equals $1$ or $3$.
In addition, it is easy to verify that under \eqn{3.5} the gap conditions on the colored parts 
of Type--$1$ partitions become identical with the gap conditions governing $A(N)$ in 
Theorem~\ref{thm:1}. Therefore, summing  over $i,j,k$ one immediately obtains 
Theorem~\ref{thm:1} from Theorem~\ref{thm:2}. 

To prove Theorem~\ref{thm:2}, Alladi, Andrews, Gordon stated it in the following form,
which they termed as a \textit{Key Identity}:
\begin{equation} 
   \sum_{\substack{i,j,k\\\text{constraints}}}
   \frac{q^{T_t+T_{ab}+T_{ac}+T_{bc-1}}
   (1-q^a+q^{a+bc})} {(q)_a(q)_b(q)_c(q)_{ab}(q)_{ac}(q)_{bc}}=
   \frac{q^{T_i+T_j+T_k}}{(q)_i(q)_j(q)_k},
\label{eq:3.6}
\end{equation}
where
$$
t=a+b+c+ab+ac+bc.
$$

Alladi and Berkovich \cite{AB} proved the following polynomial version of \eqn{3.6}
\begin{align}
   &\sum_{\substack{i,j,k\\ \text{constraints}}}
   q^{T_t+T_{ab}+T_{ac}+T_{bc-1}} \notag\\
   & \quad \times \Bigg\{ q^{bc} \qBin{L-t+a}{a}{q} \qBin{L-t+b}{b}{q}
   \qBin{L-t+c}{c}{q} \qBin{L-t}{ab}{q} \qBin{L-t}{ac}{q} \qBin{L-t}{bc}{q} \notag\\
   & \qquad \quad + \qBin{L-t+a-1}{a-1}{q} \qBin{L-t+b}{b}{q}
   \qBin{L-t+c}{c}{q} \qBin{L-t}{ab}{q} \qBin{L-t}{ac}{q} \qBin{L-t}{bc-1}{q} \Bigg\} \notag\\
   & \quad = q^{T_i+T_j+T_k}
   \qBin{L-i}{j}{q} \qBin{L-j}{k}{q} \qBin{L-k}{i}{q}.
\label{eq:3.7}
\end{align}
Actually in \cite{AB}, a doubly bounded version of \eqn{3.6} was established. Subsequently, 
the full triply bounded refinement of \eqn{3.6} was found and proven by Berkovich and Riese 
in \cite{BR}. However, for our purposes here, only the singly bounded refinement \eqn{3.7} is 
required. A partition theoretic interpretation of \eqn{3.7}, given in \cite{AB} is as follows:
\begin{thm} (Alladi, Berkovich) \\
Let $G_L(N,a,b,c,ab,ac,bc)$ be defined as $G(N,a,b,c,ab,ac,bc)$
with the additional constraint that no part exceeds $\mathbf{C}_L$.
Let $P_L(N,i,j,k)$ be defined as $P(N,i,j,k)$ with the additional constraints
\[
  \lambda(\mathbf{A})\le\mathbf{A}_{L-k},\quad
  \lambda(\mathbf{B})\le\mathbf{B}_{L-i},\quad
  \lambda(\mathbf{C})\le\mathbf{C}_{L-j},
\]
where $\lambda(\mathbf{A})$ denotes the largest part in color $\mathbf{A}$, and
$\lambda(\mathbf{B}),~\lambda(\mathbf{C})$ have the analogous interpretation.
Then
\begin{equation} 
   \sum_{\substack{i,j,k \\\text{constraints}}} G_L(N,a,b,c,ab,ac,bc)=P_L(N,i,j,k).
\label{eq:3.8}
\end{equation}
\label{thm:3}
\end{thm}

\medskip

\section{Theorem~\ref{thm:3} as a weighted partition identity}

\medskip

Considering transformations similar to \eqn{3.5}, Alladi \cite{Al1}, \cite{Al2}, and more 
recently Alladi and Berkovich \cite{AB1}, obtained weighted partition reformulations of 
G\"ollnitz's theorem. We will now follow the well-trodden path in \cite{Al1}, \cite{Al2} and
\cite{AB1}, and reformulate Theorem~\ref{thm:3} as a weighted partition identity.

To this end we will multiply both sides in \eqn{3.8} by $A^i B^j C^k q^N$ and sum over 
$i,j,k,N$ to obtain
\begin{align}
   &\sum_{\substack{a,b,\ldots,bc\ge 0 \\ N\ge 0}}
   A^{a+ab+ac} B^{b+ab+bc} C^{c+ac+bc} q^N G_L(N,a,b,c,ab,ac,bc)= \notag\\
   & \sum_{i,j,k\ge 0} q^{T_i+T_j+T_k} A^i B^j C^k
   \qBin{L-i}{j}{q} \qBin{L-j}{k}{q} \qBin{L-k}{i}{q}. 
\label{eq:4.1}
\end{align}
Actually, the combinatorial meaning of the above operation is quite natural.
Our emphasis now is to count a part occuring in color $\mathbf{A}$, or $\mathbf{B},\ldots$, or
$\mathbf{BC}$ with weight $A$, or $B,\ldots$, or $BC$, respectively.
For purposes of brevity and convenience we denote by $\mathbf{X}$ one of the six colors 
above and by $X$ its corresponding weight. For example, if $\mathbf{X}=\mathbf{AB}$, then   
$X=AB$. Here, $AB$ actually means $A$ times $B$, as do $AC$ and $BC$.

Next, we make the following crucial observation. In Type--$1$ partitions, if two parts differ 
by $\ge 2$, then colors could be assigned to these parts in any way we please. Thus, there is 
total ``independence" in assigning these colors. So, we may proceed as follows.

Consider a partition $\tilde\pi$ into distinct (uncolored) parts. Decompose $\tilde\pi$ into 
chains $\chi$, 
where a chain is a maximal run of consecutive integers. Owing to the ``independence" mentioned 
above, the weight $\omega(\tilde\pi)$ of $\tilde\pi$ may be defined multiplicatively as
\begin{equation} 
   \omega(\tilde\pi)=\prod_\chi \omega(\chi),
\label{eq:4.2}
\end{equation}
where the product is over all chains $\chi\in\tilde\pi$. In the Appendix, we will determine 
the weights of these chains as polynomials in $A,B,C$ such that the resulting weight 
$\omega(\tilde\pi)$ will be the same as the one we obtain if we attach all possible colors 
from the set $\{\mathbf{A},\ldots,\mathbf{BC}\}$ to the parts of $\tilde\pi$ in a manner 
consistent with the color-gap conditions for Type--$1$ partitions and then sum over all 
allowed color assignments counting each part in color $\mathbf{X}$ with weight $X$. 
In other words, the resulting weight will be
\begin{equation} 
   \omega(\tilde\pi)=\sum\prod_{i=1}^{l(\tilde\pi)} X_i,
\label{eq:4.3}
\end{equation}
where $l(\tilde\pi)$ is the number of parts of $\tilde\pi$, $X_i$ is the weight of $i$-th part 
of $\tilde\pi$ in a certain color assignment, and sum is over all allowed color assignments. 
For example, if $\tilde\pi_1$ 
represents partition $2+4$, then $\omega(\tilde\pi_1)=(A+B+C+A B+A C+B C)^2$, and if 
$\tilde\pi_2$ represents partition $1+2$, then $\omega(\tilde\pi_2)=A(BC+B+C)+B(B+C)+C^2$.
This way we will now establish the following
\begin{thm}
Let $D_L$ denote the set of all partitions $\tilde\pi$ into distinct parts $\le L$ with weights 
$\omega(\tilde\pi)$, as in \eqn{4.2}, \eqn{A.14} and \eqn{A.15}. Then
\begin{equation}
   \sum_{\tilde\pi\in D_L} \omega(\tilde\pi) q^{|\tilde\pi|}=
   \sum_{i,j,k\ge 0} q^{T_i+T_j+T_k} A^i B^j C^k
   \qBin{L-i}{j}{q} \qBin{L-j}{k}{q} \qBin{L-k}{i}{q}, 
\label{eq:4.4}
\end{equation}
where $|\tilde\pi|$ is the sum of parts of $\tilde\pi$.
\label{thm:4}
\end{thm}

Since
\begin{equation}
   \lim_{L\rightarrow\infty} \qBin{L}{i}{q}=\frac{1}{(q)_i},
\label{eq:4.5}
\end{equation}
we see that \eqn{4.4} together with \eqn{2.21} yields the following
\begin{crl}
Let $D$ denote the set of all partitions $\tilde\pi$ into distinct parts with weights 
as in \eqn{4.2}, \eqn{A.14} and \eqn{A.15}. Then
\begin{equation}
   \sum_{\tilde\pi\in D} \omega(\tilde\pi) q^{|\tilde\pi|}=(-Aq,-Bq,-Cq;q)_\infty.
\label{eq:4.6}
\end{equation}
\label{crl:1}
\end{crl}
As can be seen from the formulas \eqn{A.14}, \eqn{A.15} derived in the Appendix, weights 
$\omega(\chi)$, in general, are somewhat unwieldy, although there is a certain pattern as 
indicated in \eqn{A.15}. However, if one sets $C=-1$ in these formulas, it turns out that 
there is a dramatic collapse. More precisely, with $C=-1$ we have (see Appendix for derivation)
\begin{equation} 
\omega(\chi)=
\begin{cases} (-1)^{l(\chi)}(1-AB), &\mbox{ if }s(\chi)>1, \\ 
(-1)^{l(\chi)}\big\{1+\sum_{i=1}^{l(\chi)}((-A)^i+(-B)^i)\big\}, &\mbox{ if }s(\chi)=1, 
\end{cases}  
\label{eq:4.7}
\end{equation}
where $l(\chi)$ is the length (number of parts) of the chain $\chi$ and $s(\chi)$ is the 
least part of $\chi$. Furthermore, setting $A=\frac{1}{B}=z,C=-1$ in \eqn{4.7} yields 
\begin{equation} 
\omega(\chi)=
\begin{cases} \frac{z^{-l(\chi)}+z^{1+l(\chi)}}{1+z}, &\mbox{ if }s(\chi)=1, \\ 
0, &\mbox{ otherwise}. \end{cases}  
\label{eq:4.8}
\end{equation}
Now, since
\begin{equation} 
|\chi|=1+2+\ldots+l(\chi)=T_{l(\chi)}, \qquad \mbox{ if } s(\chi)=1,
\label{eq:4.9}
\end{equation}
we obtain with the aid of \eqn{4.4} and \eqn{4.8}
\begin{equation}
   \sum_{l=0}^L q^{T_l}\frac{z^{-l}+z^{1+l}}{1+z} =
   \sum_{i,j,k\ge 0} q^{T_i+T_j+T_k} z^{i-j} (-1)^k
   \qBin{L-i}{j}{q} \qBin{L-j}{k}{q} \qBin{L-k}{i}{q}, 
\label{eq:4.10}
\end{equation}
which is \eqn{1.12}. Note that with $z=1$ \eqn{4.10} becomes
\begin{equation}
   \sum_{l=0}^L q^{T_l} =
   \sum_{i,j,k\ge 0} q^{T_i+T_j+T_k}(-1)^k
   \qBin{L-i}{j}{q} \qBin{L-j}{k}{q} \qBin{L-k}{i}{q}, 
\label{eq:4.11}
\end{equation}
a polynomial version of a well-known formula of Gauss
\begin{equation}
   \sum_{l=0}^\infty q^{T_l}=(-q)_\infty^2(q)_\infty=(-q)_\infty(q^2;q^2)_\infty=
   \frac{(q^2;q^2)_\infty}{(q;q^2)_\infty}.
\label{eq:4.12}
\end{equation}
It is instructive to compare \eqn{4.11} with two formulas of Shanks \cite{S}
\begin{equation}
   \sum_{l=0}^{2L-1} q^{T_l}=\sum_{i=0}^{L-1} q^{i(2L+1)}
   \frac{(q^{2i+2};q^2)_{L-i}}{(q^{2i+1};q^2)_{L-i}},
\label{eq:4.13}
\end{equation}
\begin{equation}
   \sum_{l=0}^{2L} q^{T_l}=\sum_{i=0}^L q^{i(2L+1)}
   \frac{(q^{2i+2};q^2)_{L-i}}{(q^{2i+1};q^2)_{L-i}}.
\label{eq:4.14}
\end{equation}

There is another case that leads to dramatic simplifications in weights $\omega(\chi)$ namely,
$$
   A=-B=z, \qquad C=1.
$$
Here, we deduce from \eqn{A.14} and \eqn{A.15} that
\begin{equation} 
   \omega(\chi)=
   \begin{cases} 1, &\mbox{ if }s(\chi)=1, \\ 
   1-z^2, &\mbox{ if }s(\chi)>1. \end{cases}  
\label{eq:4.15}
\end{equation}
In addition, for $A=-B=z, C=1$ the product in \eqn{4.6} can be interpreted as
\begin{equation}
   (-zq,zq,-q;q)_\infty=(z^2q^2;q^2)_\infty (-q)_\infty=
   \frac{(z^2q^2;q^2)_\infty}{(q;q^2)_\infty}=   
   \sum_{n,j\ge 0}E(n,j) q^n(-z^2)^j,
\label{eq:4.16}
\end{equation}
where $E(n,j)$ denotes the number of partitions of $n$ with exactly $j$ distinct even parts 
(all other parts being odd). And so, combining \eqn{4.6}, \eqn{4.15}, and \eqn{4.16} we 
arrive at
\begin{equation}
   \sum_{k\ge 0}V(n,k)(1-z^2)^k=\sum_{j\ge 0}E(n,j)(-z^2)^j,
\label{eq:4.17}
\end{equation}
where $V(n,k)$ is the number of partitions of $n$ into distinct parts with exactly $k$ chains 
of consecutive integers $\neq 1$. Comparing coefficients in \eqn{4.17} we rediscover Andrews's 
generalization of Euler's partition theorem \cite{A2}
\begin{equation}
   E(n,j)=\sum_{k\ge j} \left(\begin{array}{c} k \\ j \end{array}\right) V(n,k),
\label{eq:4.18}
\end{equation}
also discussed in \cite{AG}.

We remark that \eqn{4.18} with $j=0$ becomes Euler's theorem:
\begin{equation}
   p_o(n)=p_d(n),
\label{eq:4.19}
\end{equation}
where $p_o(n)$ and $p_d(n)$ denote the number of partitions of $n$ into odd and distinct parts, 
respectively. It will be seen in the next section that $q$-version of \eqn{4.8} is
Lebesgue's identity as stated in \eqn{4.20}.

\medskip

\section{Two-parameter extension of Jacobi's triple product identity}

\medskip

In this section we will employ Corollary~\ref{crl:1} with $C=-1$ and weights given in \eqn{4.2} 
and \eqn{4.7} to derive the following generalization of Jacobi's identity \eqn{1.5}:
$$
   \sum_{i,j\ge 0} (-1)^i q^{T_{i+j}+i} \frac{(AB)_i}{(q)_i}
   \left\{\frac{A^{1+j}}{1+A}+\frac{B^{1+j}}{1+B}\right\}+
$$
\begin{equation} 
   \frac{(1-AB)}{(1+A)(1+B)}
   \sum_{i,j\ge 0} (-1)^{i+j} q^{T_{i+j}+i} \frac{(AB)_i}{(q)_i}= (-Aq,-Bq,q;q)_\infty.
\label{eq:5.1}
\end{equation}
Note that if we set $A=\frac{1}{B}=z$, then $1-AB=0$, and, therefore, we may disregard the 
second double on the LHS of \eqn{5.1}. In addition, the first double sum in \eqn{5.1} becomes 
single-fold sum because $(AB)_i=(1)_i\neq 0$ only if $i=0$. Thus, \eqn{5.1} reduces in this 
case to \eqn{1.5}. To establish \eqn{5.1} we begin by proving the following
\begin{lem} 
Let $g_i(y,q)$ denote the generating function of unrestricted partitions $\pi$ into exactly 
$i$ parts such that $\pi$ is counted with weight $(1-y)^{\nu_d(\pi)}(-1)^i$, where
$\nu_d(\pi)$ is the number of different parts of $\pi$.
Then
\begin{equation}
   g_i(y,q)=(-q)^i\frac{(y;q)_i}{(q)_i}.
\label{eq:5.2}
\end{equation}
\label{lem:1}
\end{lem}
To prove \eqn{5.2} we first consider partitions where all parts are equal to $j$, and $j$
occurs with frequency $f_j\ge 0$. Clearly, the generating function for these partitions is  
\begin{equation} 
   1+\sum_{f_j\ge 1} (1-y)(-\tilde t)^{f_j} q^{jf_j}=  
   \frac{1+y\tilde t q^j}{1+\tilde t q^j},
\label{eq:5.3}
\end{equation}
where we introduced the additional parameter $\tilde t$ to keep track of the total 
number of parts. And so,
\begin{equation} 
   \prod_{j=1}^\infty \frac{(1+y\tilde t q^j)}{(1+\tilde t q^j)}=
   \frac{(-\tilde t y q;q)_\infty}{(-\tilde t q;q)_\infty},
\label{eq:5.4}
\end{equation}
is the generating function for unrestricted partitions such that each part is counted 
with the weight $-\tilde t$ and each different part is counted with the weight $1-y$.
Obviously,
\begin{equation} 
   \big[\tilde t^i \big]\frac{(-\tilde t y q)_\infty}{(-\tilde t q)_\infty}=g_i(y,q),
\label{eq:5.5}
\end{equation}
where $\left[\tilde t^i\right]f(\tilde t)$ is the coefficient of $\tilde t^i$ in the 
expansion of $f(\tilde t)$ in powers of $\tilde t$. If we now expand \eqn{5.4} 
with the aid of Cauchy's identity \eqn{1.16}, we get
\begin{equation} 
   \frac{(-\tilde t y q)_\infty}{(-\tilde t q)_\infty}=
   \sum_{i=0}^\infty \tilde t^i (-q)^i \frac{(y)_i}{(q)_i},
\label{eq:5.6}
\end{equation}
from which \eqn{5.2} follows, as desired.

Analogously, expanding the product
\begin{equation} 
   \frac{1}{1-\tilde t} \prod_{j=1}^\infty \frac{1-\tilde t y q^j}{1-\tilde t q^j}=
   \frac{(\tilde t y q)_\infty}{(\tilde t)_\infty},
\label{eq:5.7}
\end{equation}
in powers of $\tilde t$, we can prove the following
\begin{lem} 
Let $h_i(y,q)$ denote the generating function of unrestricted partitions $\pi$ into exactly 
$i$ nonnegative parts such that $\pi$ is counted with weight $(1-y)^{\nu_d^\dagger(\pi)}$, 
where $\nu_d^\dagger(\pi)$ is the number of different positive parts of $\pi$.
Then
\begin{equation}
   h_i(y,q)=\frac{(yq)_i}{(q)_i}.
\label{eq:5.8}
\end{equation}
\label{lem:2}
\end{lem}
Let us now consider partition $\tilde\pi_t\in D^t$, where $D^t$ denotes the set of all 
partitions into distinct parts such that total number of parts equals $t$. Let us further 
assume that $\tilde\pi_t$ can be decomposed into chains of consecutive integers as
\begin{equation}
   \chi_0\cup\chi_1\cup\chi_2\cup\chi_3\ldots,
\label{eq:5.9}
\end{equation}
such that $s(\chi_0)=1,l(\chi_0)=j\ge 0$ and $s(\chi_{m>0})\neq 1$. Now if $i\ge 0$ 
denotes the total number of parts of $\chi_1\cup\chi_2\cup\chi_3\cup\ldots$, then it is 
obvious that 
\begin{equation}
   t=i+j.
\label{eq:5.10}
\end{equation}
We now subtract $1$ from the smallest part of $\tilde\pi_t$, $2$ from the second smallest 
part of $\tilde\pi_t,\ldots$, $t$ from the largest part of $\tilde\pi_t$. We call this process 
the Euler subtraction. Since the Euler subtraction can be easily reversed, it is a bijective 
process. Note that the Euler subtraction accounts for the factor $q^{T_{i+j}}$, which appears 
on the LHS of \eqn{5.1}.

Let $\tilde\pi_t^\ast$ denote the image of $\tilde\pi_t$ after the Euler subtraction, and   
$\chi_0^\ast,\chi_1^\ast,\chi_2^\ast,\ldots$ all have analogous meaning. Then it is clear that 
$\tilde\pi_t^\ast\in P_t$, where $P_t$ stands for the set of all unrestricted partitions into 
exactly $t$ nonnegative parts. Also, it is clear that $\chi_0^\ast$ represents the part $0$, 
which occurs with frequency $j=l(\chi_0)\ge 0$, and $\chi_{m>0}^\ast$ represents a positive 
part $s(\chi_m)-1-\sum_{r=0}^{m-1}l(\chi_r)$, which has frequency $l(\chi_m)$.

According to \eqn{4.7}, $\chi_0^\ast$ has the weight   
\begin{equation} 
   (-1)^j \left\{ 1+\sum_{s=0}^j (-1)^s (A^s+B^s) \right\} = (-1)^j 
   \frac{1-AB}{(1+A)(1+B)}+\frac{A^{1+j}}{1+A}+\frac{B^{1+j}}{1+B},
\label{eq:5.11}
\end{equation}
while $\chi_{m>0}^\ast$ has the weight
\begin{equation} 
   (-1)^{l(\chi_{m>0})}(1-AB). 
\label{eq:5.12}
\end{equation}
Hence, we may rewrite the LHS of \eqn{4.6} with $C=-1$ and $\omega(\tilde\pi)$ as given in 
\eqn{4.2}, \eqn{4.7} as
\begin{equation} 
   \sum_{t\ge 0}\sum_{i+j=t} q^{T_{i+j}}g_i(AB,q)
   \left\{(-1)^j \frac{1-AB}{(1+A)(1+B)}+\frac{A^{1+j}}{1+A}+\frac{B^{1+j}}{1+B}\right\}.
\label{eq:5.13}
\end{equation}
Finally, using \eqn{5.2} we derive for \eqn{4.6} with $C=-1$
$$
   \sum_{i,j\ge 0}(-1)^i q^{T_{i+j}+i}\frac{(AB)_i}{(q)_i}
   \left\{(-1)^j \frac{1-AB}{(1+A)(1+B)}+\frac{A^{1+j}}{1+A}+\frac{B^{1+j}}{1+B}\right\}=
$$
\begin{equation} 
   (-Aq,-Bq,q;q)_\infty,
\label{eq:5.14}
\end{equation}
which is essentially \eqn{5.1}, as desired.

If we repeat the above analysis for the Corollary~\ref{crl:1} with $A=-B=z,C=1$ and weights as 
in \eqn{4.2} and \eqn{4.15} we get, with the aid of Lemma~\ref{lem:2}
\begin{equation} 
   \sum_{t\ge 0} q^{T_t}h_t(z^2,q)=\sum_{t\ge 0} q^{T_t}\frac{(z^2q)_t}{(q)_t}=
   (z^2q^2;q^2)_\infty(-q)_\infty,
\label{eq:5.15}
\end{equation}
which is Lebesgue's identity, as promised at the end of \S 4.

\medskip

\section{Polynomial analogues of Lebesgue's identity}

\medskip

It is well known that Lebesgue's identity \eqn{4.20} is a special case of the $q$-Kummer sum 
\eqn{1.18}. Indeed, if in \eqn{1.18} we set $a=z^2q$ and let $b\rightarrow\infty$, we obtain 
\eqn{4.20} thanks to
\begin{equation} 
   \lim_{b\rightarrow\infty} (b;q)_n b^{-n}= (-1)^nq^{T_{n-1}}.
\label{eq:6.1}
\end{equation}
On the other hand, if in \eqn{1.18} we set $a=z^2q,b=q^{-L}$, then we derive the following 
polynomial version of Lebesgue's identity:
\begin{equation} 
   \sum_{j=0}^L \qBin{L}{j}{q} q^{T_j}(z^2 q)_j(z^2 q^{2+L+j})_{L-j}=(-q)_L(z^2q^2;q^2)_L.
\label{eq:6.2}
\end{equation}
Comparing \eqn{6.2} and \eqn{1.14} we see that these polynomial versions of \eqn{4.20} are 
radically different.

As a first step towards \eqn{1.14} we extend Lemma~\ref{lem:2} in \S 5 as follows.
\begin{lem} 
Let $h_{L,i}(y,q)$ be defined as $h_i(y,q)$ in Lemma~\ref{lem:2} with an additional constraint 
that no part exceeds $L$.
Then
\begin{equation} 
   h_{L,i}(y,q)=\sum_{j\ge 0} q^{T_j} (-y)^j \qBin{L+i-j}{i}{q} \qBin{i}{j}{q}.
\label{eq:6.3}
\end{equation}
\label{lem:3}
\end{lem}

To prove \eqn{6.3} we observe that
\begin{equation} 
   h_{L,i}(y,q)=\left[\tilde t^i\right]\frac{(qy\tilde t)_L}{(\tilde t)_{L-1}}.
\label{eq:6.4}
\end{equation}
To expand the ratio $\frac{(qy\tilde t)_L}{(\tilde t)_{L-1}}$ in powers of $\tilde t$ we can 
employ formula \eqn{1.11} along with another identity of Euler
\begin{equation} 
   \frac{1}{(\tilde t)_{L+1}}=\sum_{n=0}^\infty \tilde t^n \qBin{n+L}{n}{q}
\label{eq:6.5}
\end{equation}
to get
\begin{equation} 
   \frac{(qy\tilde t)_L}{(\tilde t)_{L+1}}=    
   \sum_{n,m\ge 0}(-y)^n \tilde t^{n+m} q^{T_n} \qBin{L}{n}{q} \qBin{m+L}{m}{q}.
\label{eq:6.6}
\end{equation}
Using \eqn{6.4} and \eqn{6.6} it is straightforward to verify \eqn{6.3}. We are now ready 
to deduce the identity \eqn{1.14} from Theorem~\ref{thm:4} with $A=-B=z,C=1$ and weights 
as given in \eqn{4.2}, \eqn{4.15}. Since our analysis will be very similar to that carried in 
\S 5, we will only sketch it here.

Let us consider partition $\tilde\pi_{L,i}\in D_{L,i}$, where $D_{L,i}$ denoted the set of all 
partitions into $i$ distinct parts, each $\le L$. As before we decompose  $\tilde\pi_{L,i}$ 
into chains of consecutive integers (none exceeding $L$) as
$$
   \tilde\pi_{L,i}=\chi_0\cup\chi_1\cup\chi_2\cup\ldots,
$$
where $s(\chi_0)=1,l(\chi_0)\ge 0,s(\chi_{m>0})>1$ and  
\begin{equation}
   i=l(\chi_0)+l(\chi_1)+l(\chi_2)+\ldots\mbox{ .}
\label{eq:6.7}
\end{equation}
We now perform the Euler subtraction, explained in \S 5. Note that this subtraction gives rise 
to the factor $q^{T_i}$ in the LHS of \eqn{1.14}. Let us introduce symbols 
$\tilde\pi_{L,i}^\ast, \chi_m^\ast$, which denote image of $\tilde\pi_{L,i}$ and $\chi_m$, 
respectively, after the Euler subtraction. It is clear that
\begin{equation}
  \tilde\pi_{L,i}^\ast\in P_{L-i,i},
\label{eq:6.8}
\end{equation}
where $P_{L,i}$ denotes the set of all unrestricted partitions into $i$ nonnegative parts, 
each $\le L$. In addition, $\chi_0^\ast$ represents part $=0$, which occurs with the 
frequency $l(\chi_0)\ge 0$, and $\chi_m^\ast>0$ represents positive part 
$s(\chi_{m>0})-1-\sum_{r=0}^{m-1}l(\chi_r)$, whose frequency is $l(\chi_{m>0})$.

Now, according to \eqn{4.15} we need to count $\chi_0^\ast$ with weight $=1$, and 
$\chi_{m>0}^\ast$ with the weight $=(1-z^2)$. Thus, we can rewrite the LHS of \eqn{4.4} as
\begin{equation} 
   \sum_{i\ge 0} q^{T_i}h_{L-i,i} (z^2,q).
\label{eq:6.9}
\end{equation}
Recalling \eqn{6.3}, we have with $A=-B=z, C=1$ and weights as given in \eqn{4.2} 
and \eqn{4.15}, the following specialization of \eqn{4.4},
$$
   \sum_{i\ge 0} q^{T_i} \sum_{j\ge 0} q^{T_j}(-z^2)^j \qBin{L-j}{i}{q} \qBin{i}{j}{q}=
$$
\begin{equation}
   \sum_{i,j,k\ge 0} (-1)^j q^{T_i+T_j+T_k} z^{i+j} 
   \qBin{L-i}{j}{q} \qBin{L-j}{k}{q} \qBin{L-k}{i}{q}, 
\label{eq:6.10}
\end{equation}
which is \eqn{1.14}, as desired.

When $L\rightarrow\infty$, \eqn{6.10} reduces to \eqn{4.20}, because \eqn{1.11} can be used 
to evaluate the inner sum on the left of \eqn{6.10}, and \eqn{2.21} can be utilized to show 
that the right hand side of \eqn{6.10} has the desired product representation.

This is not the end of the story, however. Recently, George Andrews communicated to us the 
following polynomial version of Lebesgue's identity
\begin{equation}
   \sum_{i,j\ge 0} q^{T_i+T_j}(-z^2)^j \qBin{L-j}{i}{q} \qBin{i}{j}{q}= 
   \sum_{j\ge 0} (z^2q^2;q^2)_j q^{T_{L-2j-1}} \qBin{L+1}{2j+1}{q}. 
\label{eq:6.11}
\end{equation}
Identity \eqn{6.11} can be proven by showing that both sides in \eqn{6.11} satisfy the same 
second order recurrences:
\begin{equation}
   F_L(z,q)=(1+q^L)F_{L-1}(z,q)-z^2q^LF_{L-2}(z,q),
\label{eq:6.12}
\end{equation}
and the same initial conditions
\begin{equation}
   F_0(z,q)=1, \quad F_1(z,q)=1+q.
\label{eq:6.13}
\end{equation}
Comparing \eqn{6.10} and \eqn{6.11} we infer that 
$$
   \sum_{j\ge 0} (z^2q^2;q^2)_j q^{T_{L-2j-1}} \qBin{L+1}{2j+1}{q}=    
$$   
\begin{equation}
   \sum_{i,j,k\ge 0} (-1)^j z^{i+j} q^{T_i+T_j+T_k} 
   \qBin{L-i}{j}{q} \qBin{L-j}{k}{q} \qBin{L-k}{i}{q}.  
\label{eq:6.14}
\end{equation}
It would be worthwhile to determine the precise \qhg status of \eqn{6.14}.

In a recent paper \cite{SS}, Santos and Sills discussed two polynomial analogues of 
Lebesgue's identity \eqn{4.20} with $z^2=-1,-q^{-1}$. Their results can be stated as
\begin{equation}
   \sum_{i,j\ge 0} q^{T_i+T_j} \qBin{L-j}{i}{q} \qBin{i}{j}{q}=
   \sum_{j={-\infty}}^\infty (-1)^j q^{2j^2+j}
   \mathbf{T}_1 \left(\begin{array}{c} L+1 \\ 4j+1 \end{array}\right), 
\label{eq:6.15}
\end{equation}
and
\begin{equation}
   \sum_{i,j\ge 0} q^{T_i+T_{j-1}} \qBin{L-j}{i}{q} \qBin{i}{j}{q}=
   \sum_{j={-\infty}}^\infty (-1)^j q^{2j^2}
   \mathbf{T}_1 \left(\begin{array}{c} L+1 \\ 4j+1 \end{array}\right), 
\label{eq:6.16}
\end{equation}
where $q$-trinomial coefficients $\mathbf{T}_1\left(\begin{array}{c}L\\A\end{array}\right)$ of 
Andrews and Baxter \cite{AnBx} are defined as
\begin{equation}
   \mathbf{T}_1 \left(\begin{array}{c} L \\ A \end{array}\right)=
   \sum_{\substack{r=0 \\ r\equiv L+A(mod 2)}}^{L-|A|} 
   q^{T_r} \frac{(q)_L}{(q)_{\frac{L-A-r}{2}}(q)_{\frac{L+A-r}{2}}(q)_r}.
\label{eq:6.17}
\end{equation}
We remark that \eqn{6.15}, \eqn{6.16} are special cases of generalized Gordon-G\"ollnitz 
identities found by Berkovich, McCoy, Orrick \cite{BMO} (see (2.33) with 
$\nu=2,r'=s'=0$ and (2.34) with $\nu=2,r'=0,s'=1$). It appears that polynomial identities 
\eqn{6.15}, \eqn{6.16} can not be generalized in a simple fashion to deal with all values of 
$z$ in Lebesgue's identity \eqn{4.20}. Nevertheless, it would be interesting to understand 
the relations between \eqn{6.10} with $z^2=-1,-q^{-1}$ and \eqn{6.15}, \eqn{6.16}.

\medskip

\section{Outlook}

\medskip

One of the objects of this paper has been to illustrate how the bounded Key Identity \eqn{3.7} 
introduced in \cite{AB} can be used to discover and prove unexpected polynomial versions 
of certain $q$-series identities. This account is by no means exhaustive. In our subsequent 
publications we will explore further implications of the bounded G\"ollnitz partition theorem. 
In particular, we will prove among other things the following new identity
\begin{align}
   &\sum_{n,i\ge 0} q^{\frac{3n^2-n}{2}+T_i} \qBin{L-i-2n+2}{i,n-i}{q} z^{n+i}- 
   \sum_{n,i\ge 0} q^{\frac{3n^2+n}{2}+T_i} \qBin{L-i-2n}{i,n-1-i}{q} z^{n+1+i} \notag\\
   =&\sum_{i,i,k\ge 0} q^{3(T_i+T_j+T_k)-2i-j} z^{i+j+k} 
    \qBin{L-i}{j}{q^3} \qBin{L-j}{k}{q^3} \qBin{L-k}{i}{q^3}- \notag\\ 
   zq^{3L}&\sum_{i,j,k\ge 0} q^{3(T_i+T_j+T_k)-2i-j} z^{i+j+k} 
    \qBin{L-1-i}{j}{q^3} \qBin{L-1-j}{k}{q^3} \qBin{L-1-k}{i}{q^3},
\label{eq:7.1}
\end{align}
where 
\begin{equation}
    \qBin{L}{i,j}{q} = \qBin{L}{i}{q} \qBin{L-i}{j}{q}. 
\label{eq:7.2}
\end{equation}
It is easy to verify that as $L\rightarrow\infty$ \eqn{7.1} reduces to the famous result 
of Sylvester \cite{Syl}
\begin{equation}
   \sum_{n\ge 0} q^{\frac{3n^2-n}{2}}\frac{(1+zq^{2n})(-zq)_{n-1}}{(q)_n}=(-zq)_\infty. 
\label{eq:7.3}
\end{equation}

\medskip

\section{Appendix}

\medskip

Here we derive formulas \eqn{A.14} and \eqn{A.15} for weights $\omega(\chi)$ in \eqn{4.2}, 
and formula \eqn{4.7}. \\
Let $\omega_l(\mathbf X)$ denote weight of the chain $\chi$ of length $l$ with smallest part 
$s(\chi)$ in color $\mathbf{X}\in\{\mathbf{A},\mathbf{B},\ldots,\mathbf{BC}\}$. Then using the 
color-gap conditions (see paragraph following \eqn{3.2}) and keeping in mind that each part 
in color $\mathbf{X}$ is counted with the weight $X$, it is straightforward to derive 
the following recurrences.   
\begin{align}
&\omega_l(\mathbf{C})=C\omega_{l-1}(\mathbf{C}), \label{eq:A.1} \\
&\omega_l(\mathbf{B})=B(\omega_{l-1}(\mathbf{B})+\omega_{l-1}(\mathbf{C})),\label{eq:A.2}\\
&\omega_l(\mathbf{BC})=BC(\omega_{l-1}(\mathbf{B})+\omega_{l-1}(\mathbf{C})),\label{eq:A.3}\\
&\omega_l(\mathbf{A})=A(\omega_{l-1}(\mathbf{A})+\omega_{l-1}(\mathbf{BC})+
                        \omega_{l-1}(\mathbf{B})+\omega_{l-1}(\mathbf{C})),\label{eq:A.4}\\
&\omega_l(\mathbf{AC})=AC(\omega_{l-1}(\mathbf{A})+\omega_{l-1}(\mathbf{BC})+
                          \omega_{l-1}(\mathbf{B})+\omega_{l-1}(\mathbf{C})),\label{eq:A.5}\\
&\omega_l(\mathbf{AB})=AB(\omega_{l-1}(\mathbf{AC})+\omega_{l-1}(\mathbf{A})+
                          \omega_{l-1}(\mathbf{BC})+\omega_{l-1}(\mathbf{B})+
                          \omega_{l-1}(\mathbf{C})).\label{eq:A.6}
\end{align}
For instance, to derive \eqn{A.3} we note that if we remove the smallest part in color 
$\mathbf{BC}$ from the chain $\chi$, then the new chain will have length $l-1$ and its 
smallest part will be in color $\mathbf{B}$ or $\mathbf{C}$, according to the color-gap 
conditions for the Type--$1$ partitions. Clearly, recurrences \eqn{A.1}--\eqn{A.6} 
along with the obvious initial conditions
\begin{equation}
   \omega_1(\mathbf{X})=X,\mbox{ with }\mathbf{X}\in\{\mathbf{A},\ldots,\mathbf{BC}\} 
\label{eq:A.7}
\end{equation}
determine $\omega_l(\mathbf{X})$ uniquely. It is easy to solve \eqn{A.1}--\eqn{A.6} together 
with \eqn{A.7} in sequential fashion starting with \eqn{A.1}. This way we get
\begin{align}
&\omega_l(\mathbf{C})=C^l, \label{eq:A.8} \\
&\omega_l(\mathbf{B})=\sum_{\substack{j+k=l \\ j\ge 1}} B^jC^k, \label{eq:A.9}\\
&\omega_l(\mathbf{BC})=BC\sum_{j+k=l-1} B^jC^k, \label{eq:A.10}\\
&\omega_l(\mathbf{A})=\sum_{\substack{i+j+k=l \\ i\ge 1}} 
 A^iB^jC^k+ABC\sum_{i+j+k=l-2}A^iB^jC^k, \label{eq:A.11}\\
&\omega_l(\mathbf{AC})=AC\left\{\sum_{i+j+k=l-1}A^iB^jC^k+
                             ABC\sum_{i+j+k=l-3}A^iB^jC^k\right\}+\notag\\  
                            &ABC^2\sum_{j+k=l-2}B^jC^k, \label{eq:A.12}\\
&\omega_l(\mathbf{AB})= \label{eq:A.13} \\
                 &AB\left\{\sum_{i+j+k=l-1}A^iB^jC^k + ABC\sum_{i+j+k=l-3}A^iB^jC^k+ 
                              BC\sum_{j+k=l-2}B^jC^k \right\}+ \notag\\ 
                 &A^2BC\left\{\sum_{i+j+k=l-2}A^iB^jC^k + ABC\sum_{i+j+k=l-4}A^iB^jC^k + 
BC\sum_{j+k=l-3}B^jC^k \right\}, \notag
\end{align}
where we used the convention that all summation variables are $\ge 0$, unless explicitly 
stated otherwise. 

Next, if $s(\chi)=1$, then according to the initial conditions for 
the Type--$1$ partitions the smallest part of $\chi$ may occur only in colors 
$\mathbf{A},\mathbf{B},\mathbf{C}$. Hence,
$$
   \omega(\chi)=\omega_l(\mathbf{A})+\omega_l(\mathbf{B})+\omega_l(\mathbf{C})= 
$$
\begin{equation}
   \sum_{i+j+k=l}A^iB^jC^k+ABC\sum_{i+j+k=l-2}A^iB^jC^k, \quad \mbox{if }s(\chi)=1,
\label{eq:A.14}
\end{equation}
where $l$ is the length of chain $\chi$. On the other hand, if $s(\chi)\ge 1$, then the 
smallest part of $\chi$ may occur in all six colors and, as a result, we have 
\begin{equation}
   \omega(\chi)=\omega_l(\mathbf{A})+\omega_l(\mathbf{B})+\ldots+\omega_l(\mathbf{BC})=
\label{eq:A.15}
\end{equation}
$$
   F_l(A,B,C)+A(B+C)F_{l-1}(A,B,C)+ A^2BC F_{l-2}(A,B,C), \quad \mbox{if }s(\chi)>1,
$$
where $l$ is the length of $\chi$ as before, and
\begin{equation}
   F_l(A,B,C)=\sum_{i+j+k=l}A^iB^jC^k + ABC\sum_{i+j+k=l-2}A^iB^jC^k + BC\sum_{j+k=l-1}B^jC^k. 
\label{eq:A.16}
\end{equation}

We now move on to derive \eqn{4.7}. To this end we rewrite $F_l(A,B,-1)$ as
\begin{align}
&\sum_{i+j+k=l}A^iB^j(-1)^k-\sum_{\substack{i+j+k=l, \\ i\ge 1, j\ge 1}} A^iB^j(-1)^k-
                            \sum_{\substack{j+k=l, \\ j\ge 1}} B^j(-1)^k= \notag\\ 
&\sum_{i+k=l}A^i(-1)^k+\sum_{j+k=l}B^j(-1)^k-(-1)^l-
                            \sum_{\substack{j+k=l\\j\ge 1}}B^j(-1)^k= \notag\\
&\sum_{i+k=l}A^i(-1)^k+(-1)^l-(-1)^l=(-1)^l\sum_{i=0}^l(-A)^i.
\label{eq:A.17}
\end{align}
Now, \eqn{A.15} with $C=-1$ becomes
\begin{align}
\omega(\chi)=
&(-1)^l\left\{\sum_{i=0}^l(-A)^i+(B-1)\sum_{i=1}^l(-A)^i-B\sum_{i=2}^l(-A)^i\right\}= \notag\\
&(-1)^l\left\{1+B\sum_{i=1}^l(-A)^i-B\sum_{i=2}^l(-A)^i\right\}=(-1)^l(1-AB),
\label{eq:A.18}
\end{align}
which is the first case in \eqn{4.7}. Formula \eqn{A.14} with $C=-1$ yields
\begin{align}
\omega(\chi)=&\sum_{i+j+k=l}A^iB^j(-1)^k-\sum_{\substack{i+j+k=l,\\i\ge 1,j\ge 1}}A^iB^j(-1)^k=
\notag\\ 
&\sum_{i+k=l}A^i(-1)^k+\sum_{j+k=l}B^j(-1)^k-(-1)^l= \label{eq:A.19}\\
(-1)^l&\left\{\sum_{i=0}^l(-A)^i+\sum_{j=0}^l(-B)^j-1\right\}=
(-1)^l\left\{1+\sum_{i=1}^l((-A)^i+(-B)^i)\right\},\notag
\end{align}
which is the second case in \eqn{4.7}. 

Similarly, if one sets $A=-B,C=1$ in \eqn{A.14} and \eqn{A.15}, there is a collapse as indicated 
in \eqn{4.15}. This can be shown in a manner similar to the case $C=-1$ discussed above. 

\medskip

\subsection*{Acknowledgement}

We would like to thank George E.~Andrews and Axel Riese for stimulating discussions.

\end{document}